\documentclass[graybox]{svmult}

\usepackage{type1cm}         
\usepackage{makeidx}         
\usepackage{graphicx}        
\usepackage{multicol}        
\usepackage[bottom]{footmisc}

\usepackage{newtxtext}       %
\usepackage[varvw]{newtxmath}

\makeindex

\usepackage{cleveref}
\usepackage{mathtools}
\newcommand{\update}[1]{#1}
\newcommand{\updateMath}[1]{#1}

\begin{document}

\title*{A Lattice Boltzmann Method for Non-Newtonian Blood Flow in Coiled Intracranial Aneurysms}
\titlerunning{LBM for Non-Newtonian Blood Flow in Coiled Intracranial Aneurysms}
\author{Medeea Horvat\orcidID{0009-0008-9333-1680} \and Stephan B. Lunowa\orcidID{0000-0002-5214-7245} \and Dmytro Sytnyk\orcidID{0000-0003-3065-4921} \and Barbara Wohlmuth\orcidID{0000-0001-6908-6015}}
\authorrunning{M. Horvat, S. B. Lunowa, D. Sytnyk, and B. Wohlmuth}
\institute{Medeea Horvat \and Stephan B. Lunowa \and Dmytro Sytnyk \and Barbara Wohlmuth
    \at Technical University of Munich, School of Computation, Information and Technology, Department of Mathematics, Boltzmannstraße 3, D-85748 Garching, Germany\\
    \email{medeea.horvat@tum.de}, \email{stephan.lunowa@tum.de}, \email{sytnikd@gmail.com}, \email{wohlmuth@tum.de}}
%
%
\maketitle

\abstract{
    Intracranial aneurysms are the leading cause of \update{hemorrhagic} stroke.
    One of the established  treatment approaches is the embolization induced by coil insertion.
    However, the prediction of treatment and subsequent changed flow characteristics in the aneurysm is still an open problem.
    In this work, we present an approach based on \update{a patient-specific} geometry and parameters including a coil representation as inhomogeneous porous medium.
    The model consists of the volume-averaged Navier--Stokes equations \update{for a} non-Newtonian blood rheology.
    We solve these equations using a problem-adapted lattice Boltzmann method and present a comparison between fully-resolved and volume-averaged simulations.
    The results indicate the validity of the model.
    Overall, this workflow allows for patient specific assessment of the flow due to potential treatment.
}


\section{Introduction}
Intracranial aneurysms are pathological, local deformations of the artery wall within the brain.
\update{Due to their high prevalence (about 3.2\% of the population) they are the leading cause of hemorrhagic stroke upon rupture (rupture fatality rate of 27--44\%) \cite{Pierot2013}.}
Within the last decades, the causes and treatment options have been intensively studied.
The size and location of unruptured aneurysms largely vary between patients, so that different treatment approaches were established, including e.g. contra-lateral clipping and catheter assisted coil embolization \cite{Zhao2018}.
The latter became a common treatment technique with many variations discussed in \cite{Eddleman2013}.
Nevertheless, the prediction of the treatment, including optimal coiling (occluded volume, spatial distribution) and the subsequent reduction of blood flow and wall shear stress in the aneurysm, is still an open problem, leaving the success of intervention to the medical \update{practitioner's} previous experience.

Numerical prediction of the blood-flow characteristics in coiled aneurysms requires the incorporation of patient specific geometry and parameters as well as the representation of the inhomogeneous coil distribution within the aneurysm.
For the latter, coil deployment models based on a discrete elastic rod approach were employed in \cite{Otani2017} to evaluate the effect of the coil configuration on the flow in the aneurysm.
The authors reported a significant reduction of flow momentum and kinetic energy at packing densities used in clinical practice (20--25\%).
To avoid the detailed discretization of the fluid domain with the complex coil geometry, the coil can be represented as an inhomogeneous porous medium, as done in \cite{Yadollahi2019} for straight channels and 3 patient specific geometries with computer generated coils, while the coil porosity and permeability was reconstructed from 3D X-ray images in \cite{Romero2023}.
In both approaches, the heterogeneity of the coil distribution was found to be crucial for correct flow predictions.

Beside the classical finite element and finite volume methods for flow simulations, lattice Boltzmann methods (LBM) are of major interest due to their straightforward parallelization, especially on GPUs, and the direct applicability on voxelized geometries common in medical imaging.
These methods are based on a mesoscopic description of statistical physics, i.e., the discretization of the Boltzmann equation in physical and velocity space.
For a detailed introduction to LBM, we refer to \cite{Krueger2017,Succi2018}.
Various extensions exist for different applications, in particular for non-Newtonian (blood) flow \cite{Ouared2005}, where the Casson and Carreau-Yasuda models are compared e.g. in \cite{Boyd2007,Ashrafizaadeh2009}, but also for the volume-averaged Navier--Stokes equations (VANSE) \cite{Hoecker2018,Maier2021}, which allow the fluid domain to (partially) contain a porous medium, see also \cite{He2019} for a review of earlier approaches in the context of heat transfer.
However, applications of LBM to \update{intracranial} aneurysms remain scarce \cite{Ouared2005,Hosseini2022} and limited to the fully resolved geometry.

In this work, we propose to combine the different approaches.
We use the LBM for VANSE together with a variable relaxation rate for non-Newtonian flow, such that the coil can be represented by a porous medium.
To this end, we introduce the model in \cref{sec:model} and the resulting LBM in \cref{sec:method}.
Finally, we present in \cref{sec:simulation} simulation results in a patient-specific geometry and compare the fully resolved geometry to the averaging approach for different (simulated) coiling treatment.

\section{Volume-Averaged Navier-Stokes Equations for Blood Flow}
\label{sec:model}

Let $T > 0$ be the final time and $\Omega \subset \mathbb{R}^3$ be the fluid domain.
We denote by $\Omega_T := (0,T) \times \Omega$ the time-space cylinder.
The VANSE for non-Newtonian flow are given by (see \cite[Chp. 1]{Nield2013} and references therein)
\begin{align}
    \partial_t \big(\phi \varrho\big) + \nabla\cdot\big(\phi \varrho \vec{u}\big) & = 0 & \text{in}\ \Omega_T \;, \\
    \partial_t \big(\phi \varrho \vec{u}) + \nabla \cdot \big(\phi \varrho \vec{u} \otimes \vec{u}\big)
    - \nabla \cdot \big(\updateMath{\mu(\|\dot{\vec{\gamma}}\|_F) \dot{\vec{\gamma}}}\big) + \phi \nabla p & = \varrho \vec{f}(\vec{u}) & \text{in}\ \Omega_T \;,
\end{align}
where the unknowns $\varrho : \Omega_T \to (0, \infty)$, $\vec{u} : \Omega_T \to \mathbb{R}^3$ and $p : \Omega_T \to \mathbb{R}$ denote the volume-averaged density, velocity, and pressure of the fluid, respectively.
The parameter $\phi : \Omega \to (0, 1]$ is the porosity, i.e., the volume fraction filled by fluid, while $\mu : \updateMath{\mathbb{R}} \to (0,\infty)$ is the dynamic viscosity, which depends on the \update{shear-strain rate $\dot{\vec{\gamma}} = \nabla(\phi\vec{u}) + \nabla(\phi\vec{u})^\top$
via the Frobenius norm $\|\cdot\|_F$.}
For blood, typical viscosity models are the ones by Casson and Carreau--Yasuda, see \cite{Sequeira2018} and references therein.
Here, we employ the more general Carreau--Yasuda model, viz.
\begin{align}
    \label{eq:carreau}
    \mu(\updateMath{\dot{\gamma}}) & := \mu_\infty + (\mu_0 - \mu_\infty) \big(1 + (\lambda \dot{\gamma})^a\big)^{(n-1)/a} \;,
\end{align}
where $\mu_0, \mu_\infty$ denote the limiting viscosities, $\lambda$ is the viscosity relaxation time, while $n < 1$ and $a$ are the power law index and the transition parameter, respectively.

Finally, the volume-averaged force $\vec{f}(\vec{u})$ acting on the fluid contains the interaction forces between solid and fluid,
and is typically given as
\begin{align*}
    \vec{f}(\vec{u}) = - \phi^2 \nu \tens{K}^{-1} \vec{u} - \phi^3 C_F \tens{K}^{-1/2} \|\vec{u}\| \vec{u} \;,
\end{align*}
where the first and second term model the Darcy and Forchheimer drag \update{and $\|\cdot\|$ is the Euclidean norm}.
The parameter $\nu$ denotes the kinematic viscosity in the porous medium, while the permeability $\tens{K}$ and the Forchheimer constant $C_F$ are usually defined using the (isotropic) Kozeny--Carman relation and Ergun's experiment \cite{Ergun1952}
\begin{align*}
    \tens{K} & = k \tens{Id} \;,                           &
    k        & = \frac{\phi^3 d_p^2}{150 (1 - \phi)^2} \;, &
    C_F      & = \frac{1.75}{\sqrt{150 \phi^3}} \;.
\end{align*}
Note that the VANSE coincide with the classical Navier-Stokes equations when $\phi = 1$ (pure fluid without solid fraction).

\section{Lattice Boltzmann Method for Volume-Averaged Navier-Stokes}
\label{sec:method}

We combine the D3Q27 LBM for VANSE developed in \cite{Hoecker2018,Maier2021} with a variable relaxation-rate for the non-Newtonian viscosity as discussed in \cite{Ouared2005,Boyd2007,Ashrafizaadeh2009}.
We discretize the time interval $(0,T)$ into uniform steps of size $\Delta t$ and the domain $\Omega$ into uniform cubical cells of length $\Delta x$.
\update{For D3Q27, one uses 27 discrete velocity directions $\vec{c}_i$ leading to a lattice speed of sound $c_s = 1/\sqrt{3}$ and respective lattice weights $w_i$, whose values are given e.g. in \cite{Krueger2017,Succi2018}.
The respective microscopic particle densities $f_i$} are related to the macroscopic quantities \update{(density, velocity, pressure)} by
\begin{align*}
    \varrho & := \phi^{-1} \sum_i f_i \;, &
    \vec{u} & := (\phi \varrho)^{-1} \Big(\sum_i f_i \vec{c}_i + \frac{\Delta t}{2} \varrho \vec{f}(\vec{u})\Big) \;, &
    p &:= c_s^2 \varrho \;.
\end{align*}
Since the force $\vec{f}$ has a quadratic dependence on the velocity $\vec{u}$, the equation for $\vec{u}$ can be solved for the latter.
Using the preliminary velocity $\vec{v} := (\phi \varrho)^{-1} \sum_i f_i\vec{c}_i$, the velocity $\vec{u}$ is explicitly given by
\begin{align*}
    \vec{u} &= \frac{\vec{v}}{C_0 + \sqrt{C_0^2 + C_1 \|\vec{v}\|}} \;, &
    C_0 &:= \frac{1}{2} + \frac{\phi \nu \Delta t}{4 k} \;, &
    C_1 &:= \frac{\Delta t \phi^2 C_F}{2 \sqrt{k}} \;.
\end{align*}
\begin{remark}[Velocity solution for anisotropic permeability]
    In case of \update{an} anisotropic permeability, the velocity is given by a system of interdependent quadratic equations which cannot be solved explicitly.
    Instead, one can use the fixed-point iteration \update{
    \begin{align*}
        \vec{u}_{m+1} &= \tens{X}_m^{-1}\vec{v} \;,&
        \tens{X}_m &:= \tens{A} + \|\vec{u}_{m}\| \tens{B} \;,&
        \tens{A} &:= \tens{Id} + \tfrac{\phi \nu \Delta t}{2} \tens{K}^{-1} \;,&
        \tens{B} &:=\tfrac{\phi^2 C_F \Delta t}{2}\tens{K}^{-\frac{1}{2}} \;,
    \end{align*}
    for $m\geq0$, starting e.g. with $\vec{u}_0 = \vec{u}|_{t-\Delta t}$ or $\vec{u}_0 = \vec{v}$.
    Since $\tens{X}_m$ is a sum of positive definite matrices, it is invertible and the spectral norm $\|\tens{X}_m^{-1}\|$ is bounded by $1$, so that $\|\vec{u}_{m+1}\| \leq \|\vec{v}\|$.
    Convergence follows by subtracting consecutive iterations
    \begin{align*}
        0 &= \tens{X}_m \vec{u}_{m+1} - \tens{X}_{m-1} \vec{u}_m  = \tens{X}_m (\vec{u}_{m+1} - \vec{u}_m) + (\|\vec{u}_m\| - \|\vec{u}_{m-1}\|) \tens{B} \vec{u}_m
    \shortintertext{and taking norms}
        \|\vec{u}_{m+1} - \vec{u}_{m}\|
            &= \big|\|\vec{u}_{m-1}\| - \|\vec{u}_{m}\|\big|\|\tens{X}_m^{-1}\tens{B}\vec{u}_m\|
            \leq \big\| (\tens{A} + \|\vec{v}\|\tens{B})^{-1} \|\vec{v}\|  \tens{B}\big\| \|\vec{u}_{m} - \vec{u}_{m-1}\| \;,
    \end{align*}
    since the spectral norm on the right-hand side is strictly smaller than $1$.}
    However, this iteration significantly increases the simulation time, as the repeated computation of the update is necessary in all cells of the discretization.
\end{remark}

In each time step, the particle densities $f_i$ are first collided and then streamed toward the neighboring cells.
In summary, these operations read
\begin{align*}
    f_i(t + \Delta t, \vec{x} + \Delta t \vec{c}_i)
    & = f_i(t, \vec{x}) + \Omega_i^\phi(t, \vec{x}) + \updateMath{\Delta t} \Omega_i^{\vec{f}}(t, \vec{x}) + \Omega_i^\text{SRT}(t, \vec{x}) \;.
\end{align*}
The first collision term $\Omega_i^\phi$ accounts for changes in the equilibrium due to the spatially varying porosity \cite{Hoecker2018,Maier2021},
\update{$\Omega_i^{\vec{f}}$ includes forces via the Guo forcing scheme \cite{Guo2002a}, and $\Omega_i^\text{SRT}$ is the classical single-relaxation-time collision scheme \cite[Sec.~3.2]{Krueger2017}, given by
\begin{align*}
    \Omega_i^\phi(t, \vec{x}) &:= w_i \varrho(t,\vec{x}) \big(\phi(\vec{x} + \Delta t \vec{c}_i) - \phi(\vec{x})\big) \;,\\
    \Omega_i^{\vec{f}}(t,\vec{x}) & := w_i \left(1 - \frac{\omega(t, \vec{x})}{2}\right) \left( \frac{\vec{c}_i - \vec{u}(t,\vec{x})}{c_s^2}  + \frac{(\vec{u}(t,\vec{x}) \cdot \vec{c}_i) \vec{c}_i}{c_s^4}\right) \cdot \varrho(t, \vec{x}) \vec{f}(\vec{u}(t,\vec{x})) \;, \\
    \Omega_i^\text{SRT}(t, \vec{x}) &:= -\omega(t, \vec{x}) \big(f_i(t,\vec{x}) - f_i^\text{eq}(t,\vec{x})\big) \;, \\
    f_i^\text{eq}(t,\vec{x}) & := w_i \phi(\vec{x}) \varrho(t, \vec{x}) \left( 1 + \frac{\vec{c}_i \cdot \vec{u}(t,\vec{x})}{c_s^2} + \frac{(\vec{c}_i \cdot \vec{u}(t,\vec{x}))^2}{2 c_s^4} - \frac{\vec{u}(t,\vec{x}) \cdot \vec{u}(t,\vec{x})}{2 c_s^2}
    \right) \;,
\end{align*}
where $f_i^\text{eq}$ denotes the equilibrium particle density.}
The relaxation rate \update{$\omega(t, \vec{x})$} is updated based on its relation to the dynamic viscosity by $\omega = 2 \varrho c_s^2 \Delta t (2\mu + \varrho c_s^2 \Delta t)^{-1}$.
For this, the new dynamic viscosity is computed from the Carreau-Yasuda model \cref{eq:carreau}, where the shear-strain rate \update{\renewcommand{\tens}[1]{\mathrm{#1}}%
is approximated semi-explicitly by $\dot{\gamma}(t,\vec{x}) = \omega(t - \Delta t, \vec{x}) \|\tens{\Pi}(t, \vec{x})\|_F / (c_s^2 \rho(t, \vec{x}) )$ following \cite{Ouared2005}.
Here, $\tens{\Pi}$ denotes the corrected second moment of the particle densities $f_i$ given by \cite[Sec.~6.2]{Krueger2017}
\begin{align*}
    \tens{\Pi} := \sum_i \big(f_i - f_i^\text{eq}\big) \vec{c}_i \otimes \vec{c}_i + \frac{\Delta t}{2} \varrho \big(\vec{f}(\vec{u}) \otimes \vec{u} + \vec{u} \otimes \vec{f}(\vec{u})\big) \;,
\end{align*}
which also yields the deviatoric stress tensor $\tens{\sigma}(t, \vec{x}) = (\omega(t, \vec{x}) / 2 - 1\big) \tens{\Pi}(t, \vec{x})$.}

\section{Numerical Simulations in a Realistic Aneurysm}
\label{sec:simulation}

\begin{figure}[tb]
    \begin{minipage}{\textwidth}
        \leftfigure{\includegraphics[width=54mm]{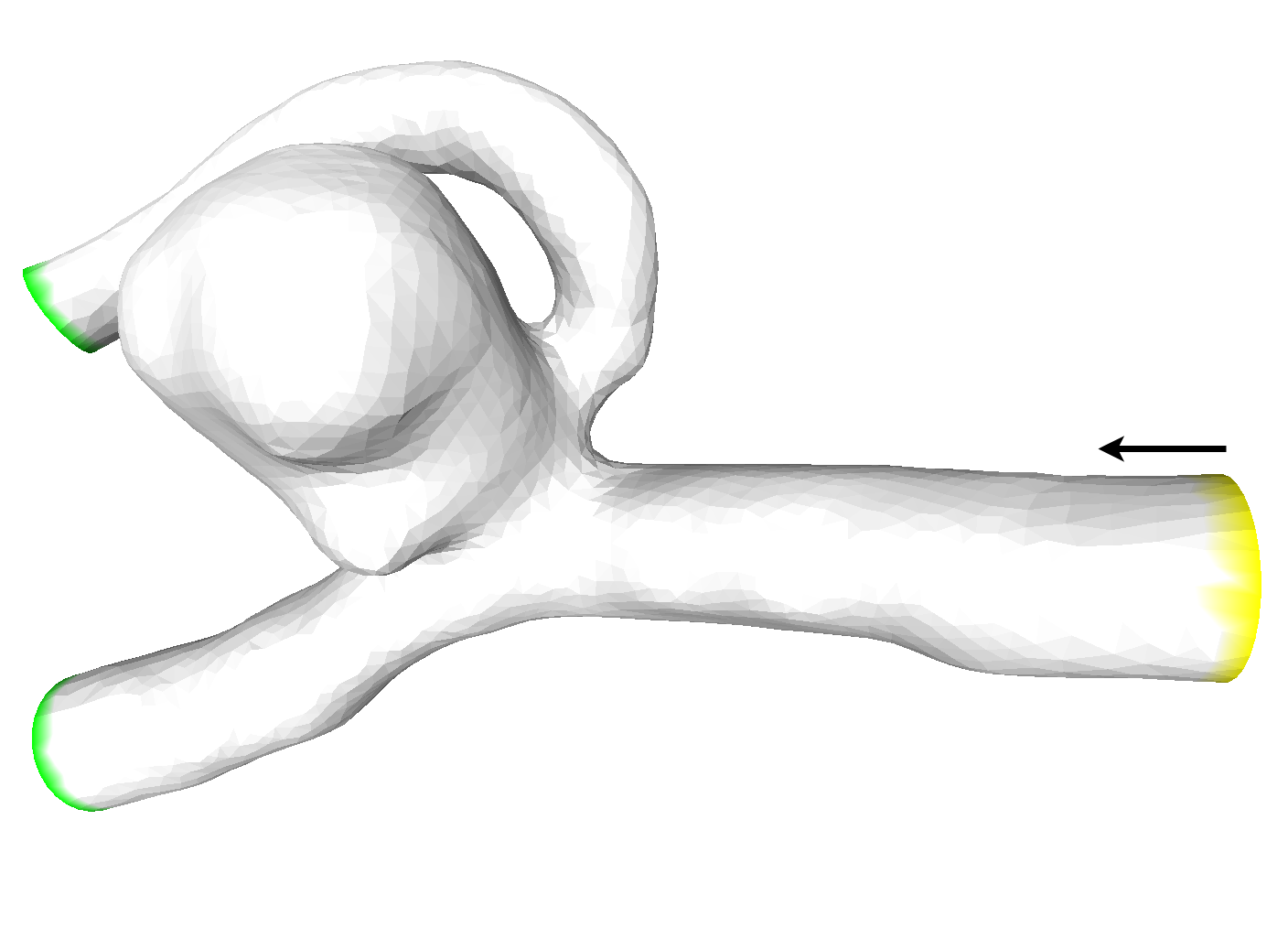}}%
        \hfill%
        \rightfigure{\includegraphics[width=54mm]{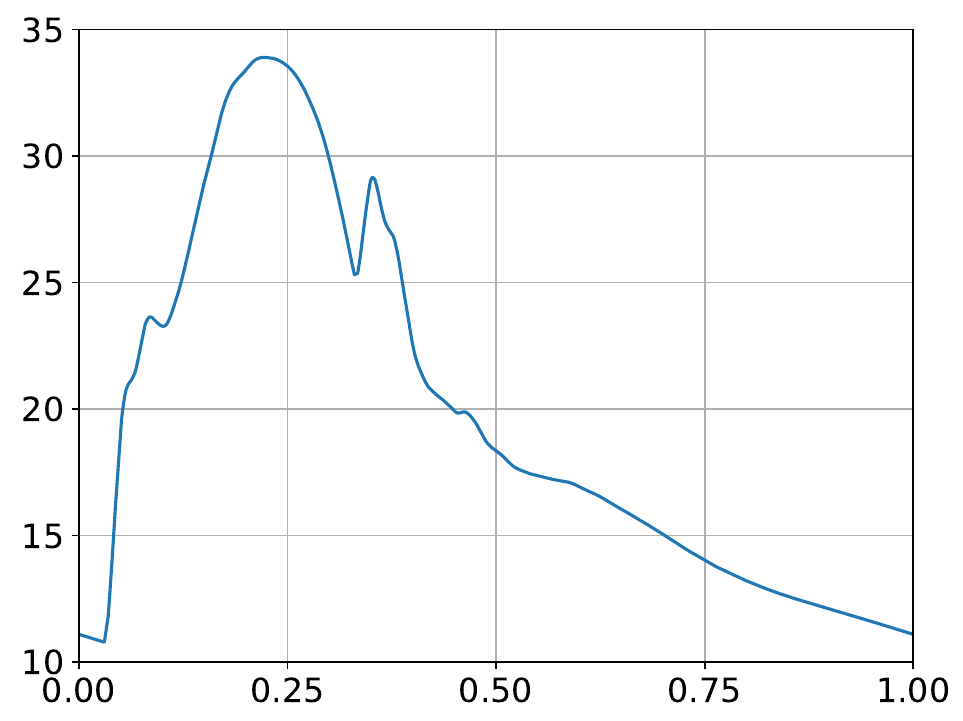}}%
    \end{minipage}
    \leftcaption{\update{Aneurysm geometry with in- and outflow marked in yellow and green, respectively.}}\label{fig:geometry}
    \rightcaption{Average inflow velocity $v$ [cm/s] over time $t$ [s] computed for arteriole 22 in \cite{Fritz2022}.}
    \label{fig:inflow}
\end{figure}

The aneurysm geometry used for the simulations is depicted in \cref{fig:geometry} and is based on an anonymized CT scan provided by the Neuro-Kopf-Zentrum of the Klinikum rechts der Isar (TU Munich).
The geometry was post-processed to enhance the mesh quality by removing scanning artefacts and by extending the in- and outflow regions to cylindrical shape to impose reasonable boundary conditions.

\begin{figure}[tb]
    \includegraphics[width=.32\textwidth, viewport=100 320 950 1110, clip=true]{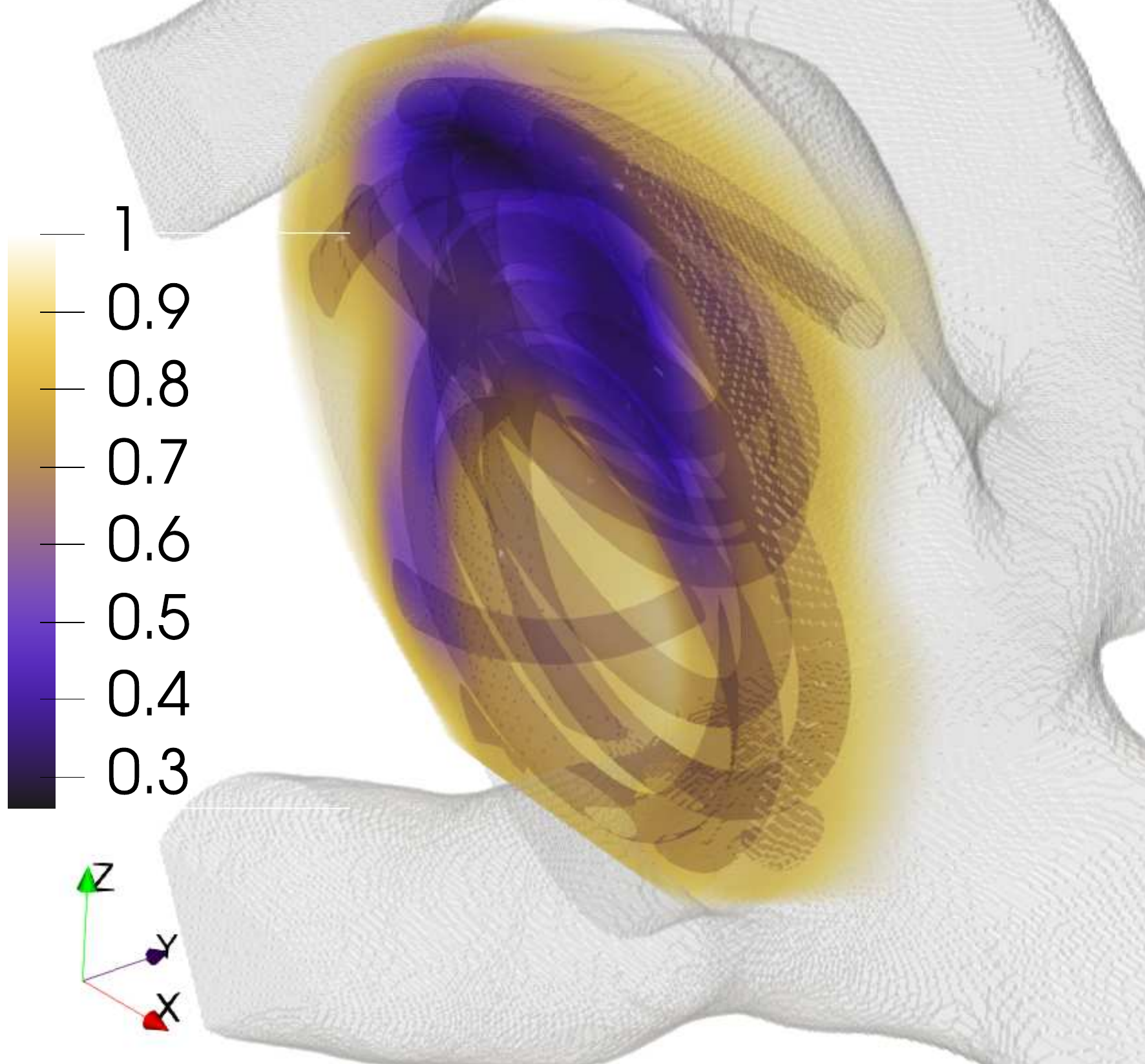}%
    \hfill%
    \includegraphics[width=.32\textwidth, viewport=100 320 950 1110, clip=true]{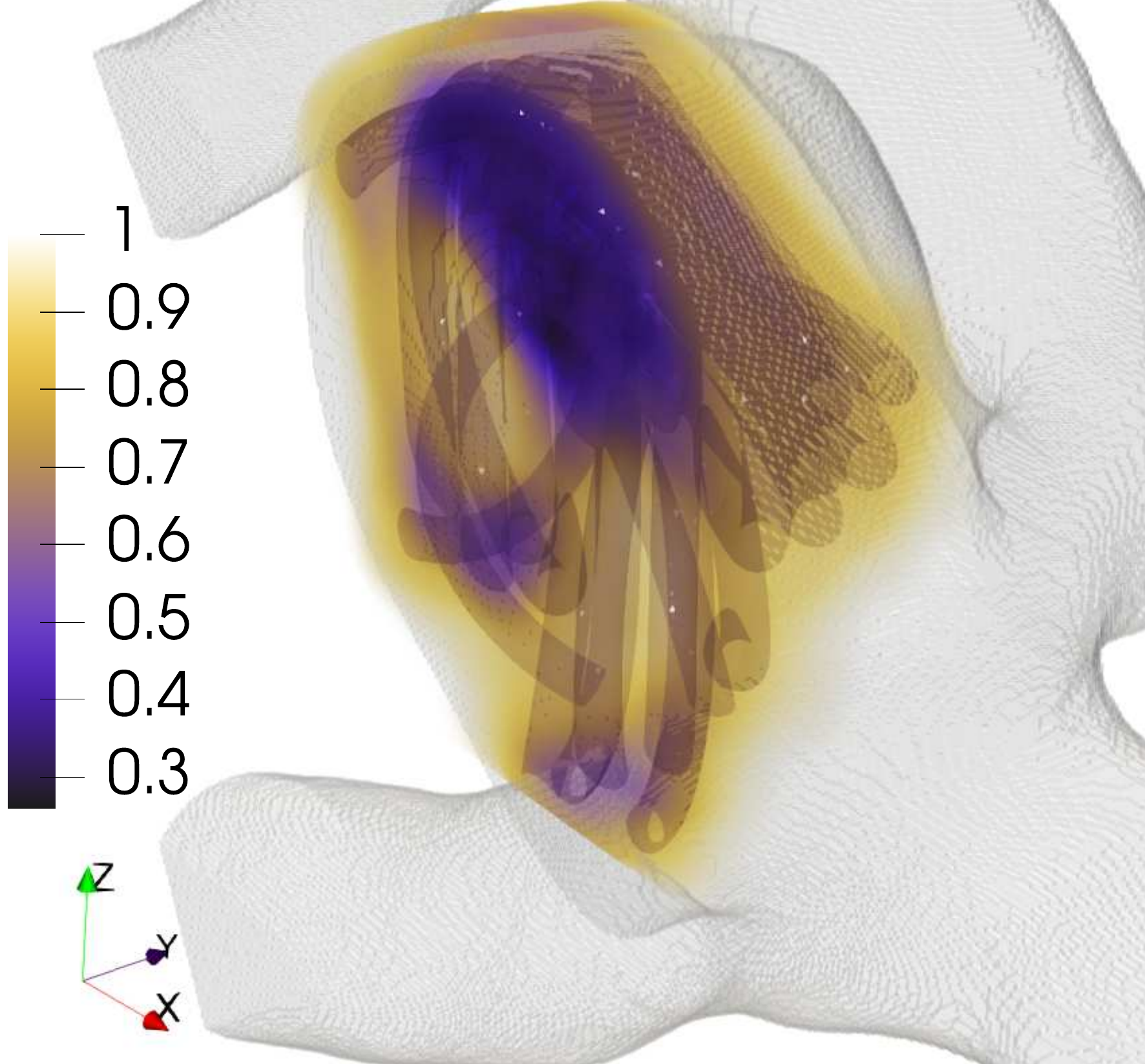}%
    \hfill%
    \includegraphics[width=.32\textwidth, viewport=100 320 950 1110, clip=true]{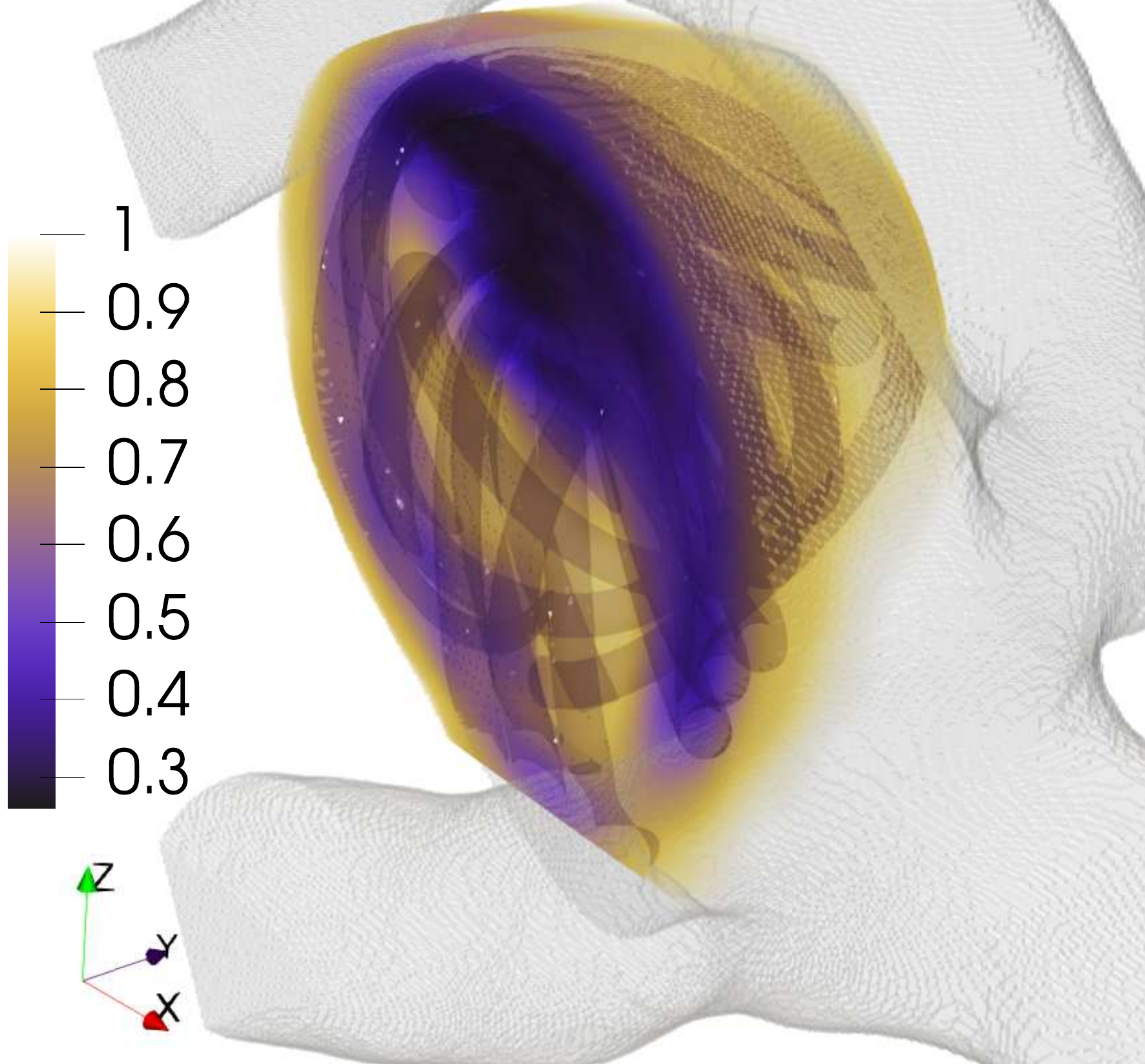}%
    \caption{Cross-section of the resulting porosity field (color) for the three coils (3D shadow) with packing density of 15\%, 20\% and 25\% (left to right).}
    \label{fig:coiling}
\end{figure}

Based on this geometry, the coiling procedure is simulated using a discrete elastic rod model discussed in \cite{Holzberger2024}.
This enables us to compare the effect of different packing densities (15\%, 20\%, 25\%), which are typically employed in surgery.
\update{T}he porosity field is obtained by windowed convolution of the voxelized coil representation with a uniform averaging kernel \update{(window size of 1.23\;mm)}.
The resulting coils and corresponding porosity fields are illustrated in \cref{fig:coiling}.

The flow simulations require physical boundary conditions given in the following.
At the arterial wall, we prescribe no-slip boundary conditions via the half-way bounce back method, see \cite[Sec.~5.3.3]{Krueger2017}.
The same holds at the coil boundary for fully resolved simulations.
At the inflow, a given velocity profile is imposed using the Zou--He method \cite{Zou1997}, while the outflow boundary condition is based on an extrapolation method \cite{Geier2015} to avoid spurious reflections.
The inflow velocity profile is based on the solution of the 1D-0D network model of arterial blood flow from \cite{Fritz2022}, which is capable of providing a realistic time-resolved mean velocity profile $v$ depicted in \cref{fig:inflow}.
Since the \update{inflow area is a planar circle, the inflow profile can be chosen} zero except in the component normal to the inflow plane, which is given by the classical Hagen--Poiseuille profile $u_{\perp}(t,r) = 2 v(t) (1 - R^{-2} r^2)$, where $r$ is the distance from the inflow center and $R$ is the inflow radius.

\begin{table}[tb]
    \renewcommand{\arraystretch}{1.25}
    \caption{Parameters for the LBM simulations}
    \label{tab:parameters}
    \begin{tabular}{l@{}crl}
        \hline\noalign{\smallskip}
        Parameter               & Symbol                & Value & Unit     \\
        \noalign{\smallskip}\svhline\noalign{\smallskip}
        Reference Density       & $\varrho_0$           & 1060  & kg/m$^3$ \\
        Inflow radius           & $R$                   & 1.649 & mm       \\
        \update{Coil wire diameter} & \update{---}      & \update{305} & \update{$\umu$m}  \\
        Final time              & $T$                   & 2.00  & s        \\
        Grid size               & $\Delta x$            & 43.9\update{5} & $\umu$m  \\
        Time-step size          & $\Delta t$            & 5.42  & $\umu$s  \\
        \noalign{\smallskip}\hline
    \end{tabular}
    \hfill
    \begin{tabular}{l@{}crl}
        \hline\noalign{\smallskip}
        Parameter               & Symbol                & Value        & Unit   \\
        \noalign{\smallskip}\svhline\noalign{\smallskip}
        \update{Porosity avg. window} & \update{---}    & \update{1.23} & \update{mm} \\
        Permeability scaling    & $d_P$                 &          600 & $\umu$m     \\
        Limit viscosities       & $\mu_0$, $\mu_\infty$ & 0.0035, 0.16 & kg/m/s \\
        Visc. relaxation time   & $\lambda$             &          8.2 & s      \\
        Visc. power-law index   & $n$                   &       0.2128 & ---    \\
        Visc. transition param. & $a$                   &         0.64 & ---    \\
        \noalign{\smallskip}\hline
    \end{tabular}
\end{table}

All simulations \update{use} the parameters listed in \cref{tab:parameters} \update{and} $\nu = \updateMath{\varrho}_0^{-1} \mu_0$ in the VANSE model.
In a pre-run phase of 10,000 time steps, starting from zero initial velocity, the inflow velocity is \update{linearly} increased to the actual value at $t = 0$\;s to avoid the formation of acoustic shock waves and thereby introduced instabilities.
Then, the actual simulations are started from the pre-run state and run for 2 heart beats with periodically repeating inflow profile to reach a close-to periodic solution.
All following reported results are at time $t = 1.23161$\;s during systole.

The implementation is based on the XLB framework \cite{Ataei2023} allowing parallel simulations on CPUs and GPUs.
The simulations were run on a workstation with AMD EPYC 7713 CPU and Nvidia RTX A6000 GPU using the JAX CUDA runtime mode (main computations were performed on GPU).
\update{This results in run times of about 11.1 and 15.5 hours for the fully-resolved and volume-averaged ones having about 5.2 million cells within the aneurysm geometry (ca.\ 140 million degrees of freedom).}

\runinhead{Fully Resolved vs. Volume-Averaged Flow.}

\begin{figure}[p]
    \begin{minipage}{\textwidth}
        \leftfigure{\includegraphics[width=.47\textwidth, viewport=370 350 1350 970, clip=true]{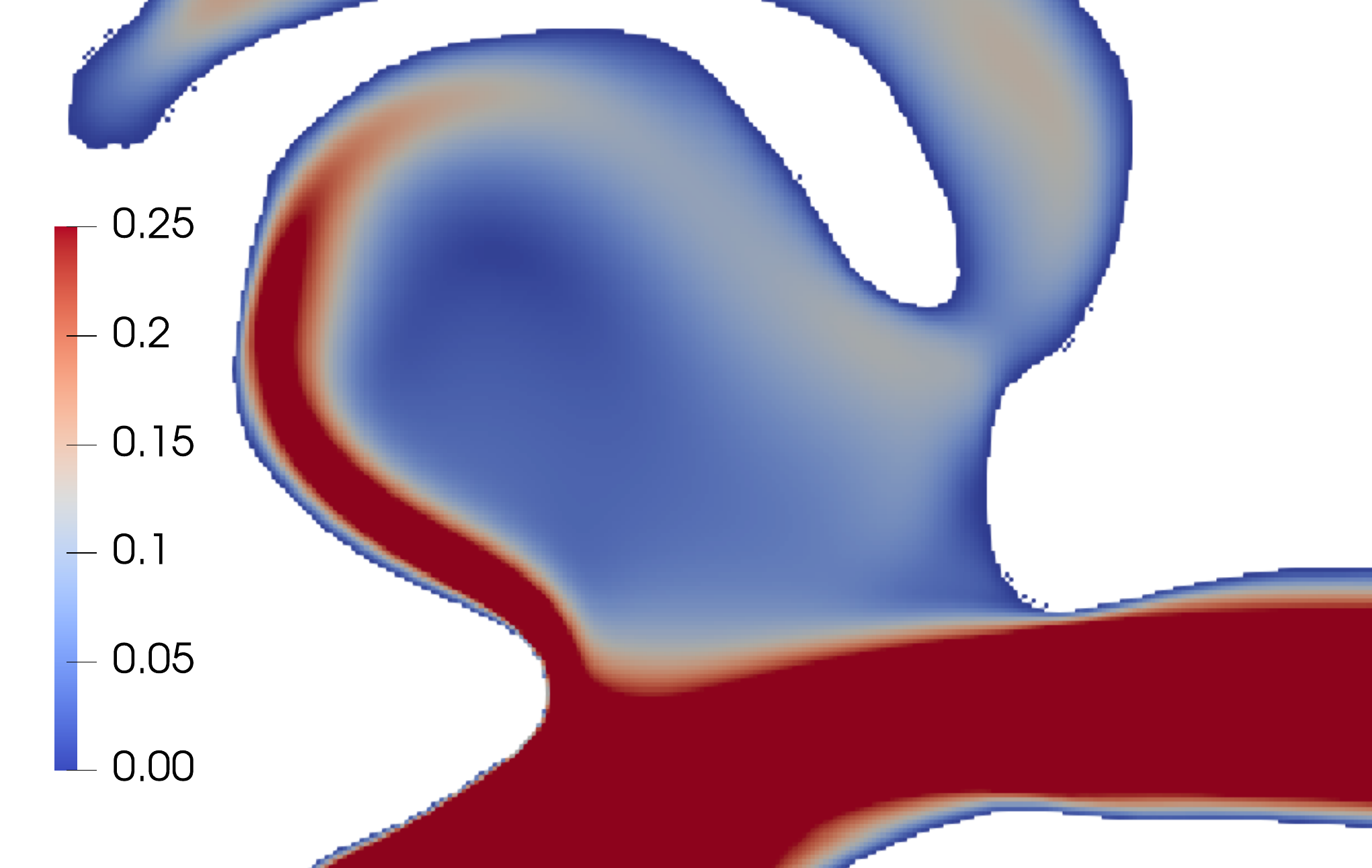}}%
        \hfill%
        \rightfigure{\includegraphics[width=.42\textwidth, viewport=0 100 1300 1020, clip=true]{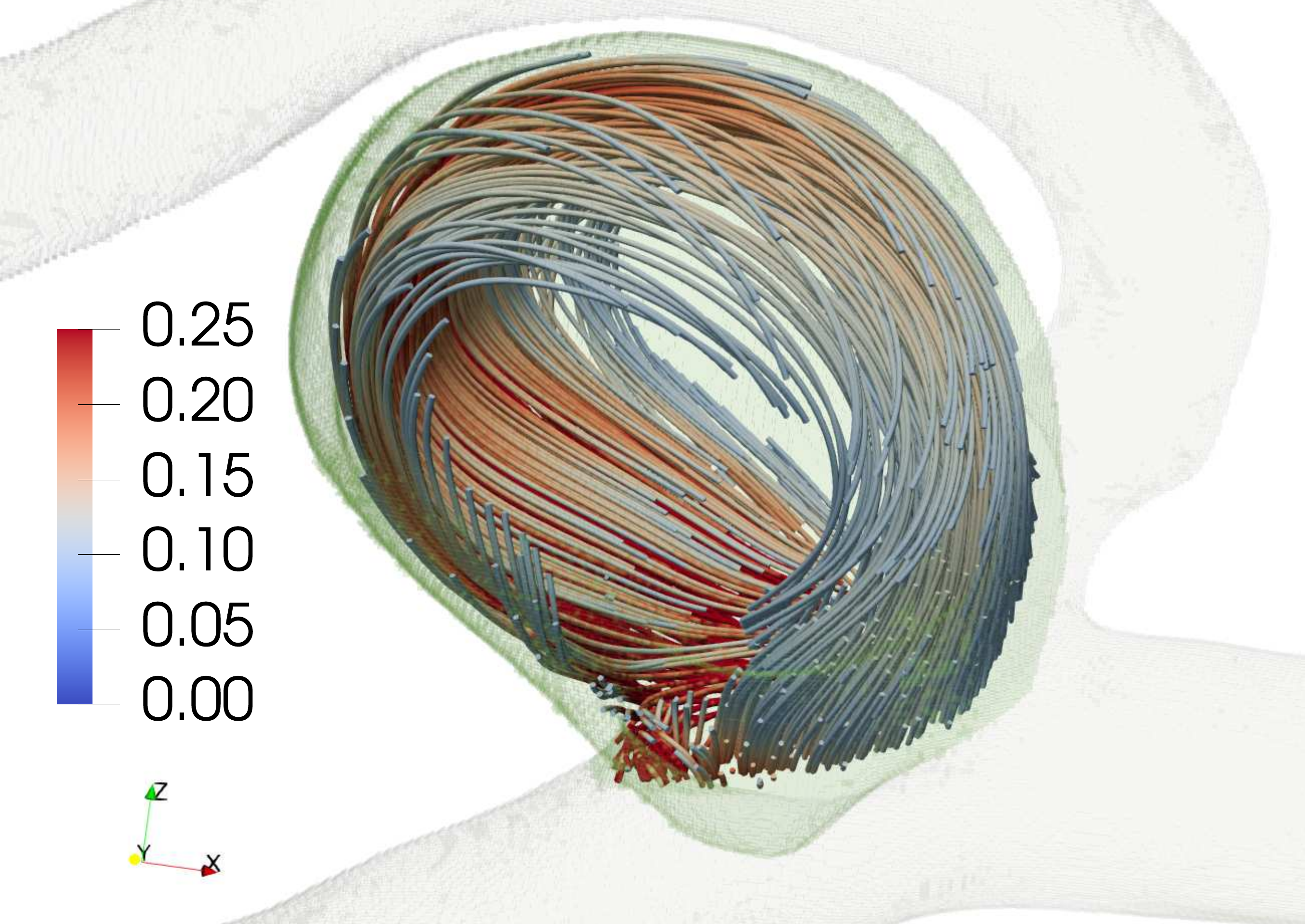}}%
    \end{minipage}
    \caption{Velocity magnitude $\|\vec{u}\|$ [m/s] during systole ($t \approx 1.23$\;s) in the aneurysm without coil plotted over a cross-section \emph{(left)} and stream lines in the cut off aneurysm region \emph{(right)}.}
    \label{fig:none}
\end{figure}
\begin{figure}[p]
    \begin{minipage}{\textwidth}
        \leftfigure{\includegraphics[width=.47\textwidth, viewport=370 350 1350 970, clip=true]{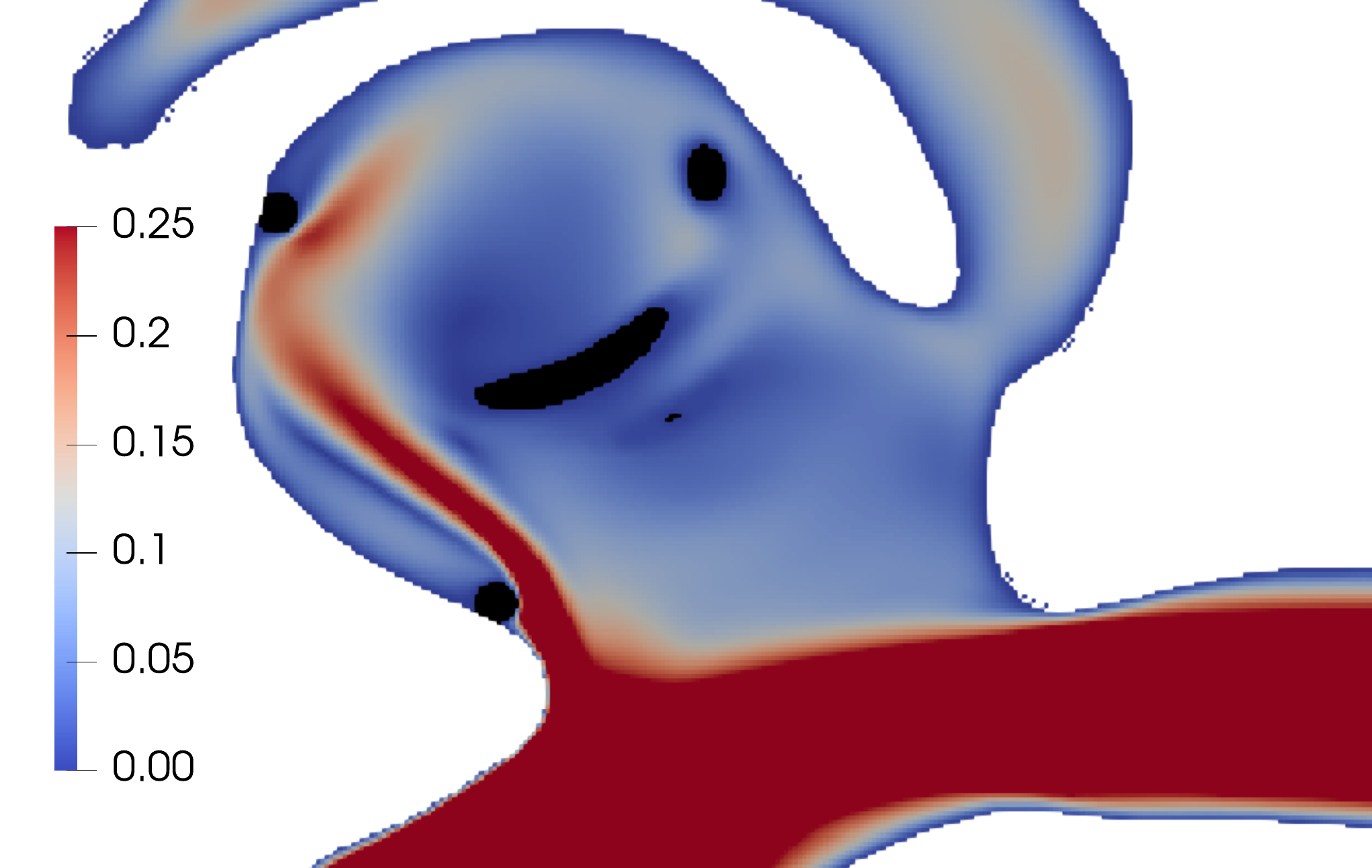}}%
        \hfill%
        \rightfigure{\includegraphics[width=.47\textwidth, viewport=370 350 1350 970, clip=true]{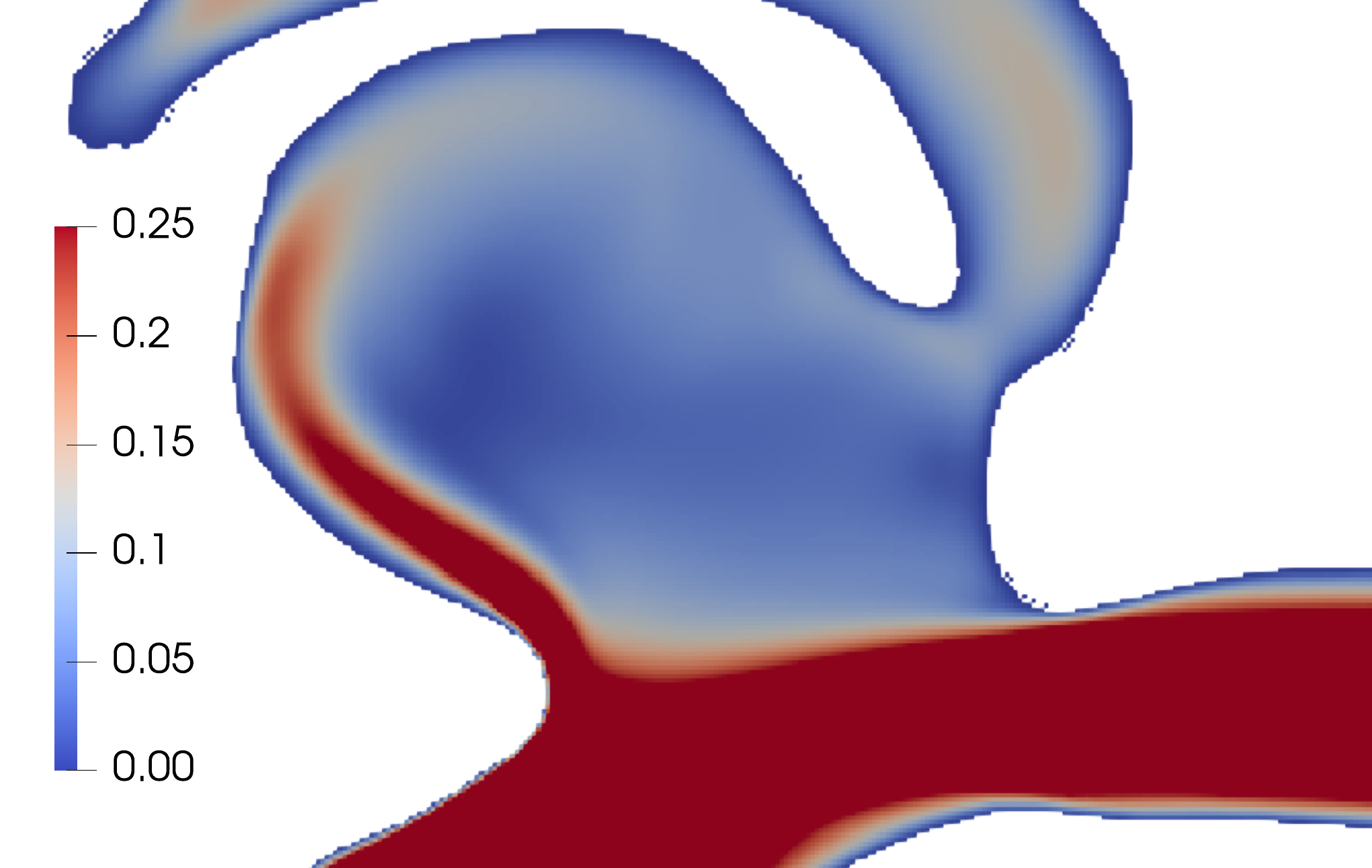}}%
    \end{minipage}
    \caption{Velocity magnitude $\|\vec{u}\|$ [m/s] during systole ($t \approx 1.23$\;s) in the aneurysm with fully resolved coil of 15\% packing density \emph{(left)} and with averaged porosity field \emph{(right)}.}
    \label{fig:vanse15}
\end{figure}
\begin{figure}[p]
    \begin{minipage}{\textwidth}
        \leftfigure{\includegraphics[width=.47\textwidth, viewport=370 350 1350 970, clip=true]{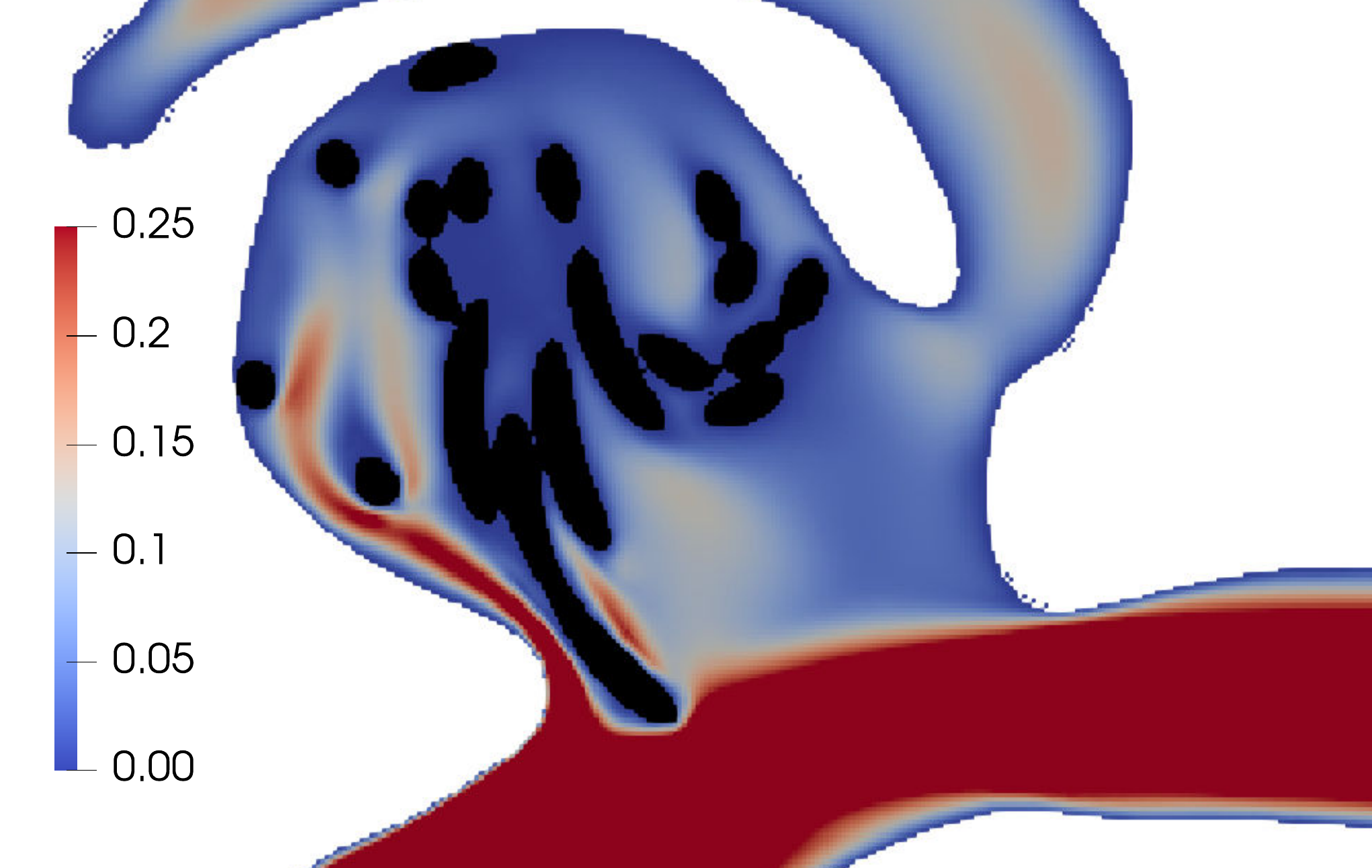}}%
        \hfill%
        \rightfigure{\includegraphics[width=.47\textwidth, viewport=370 350 1350 970, clip=true]{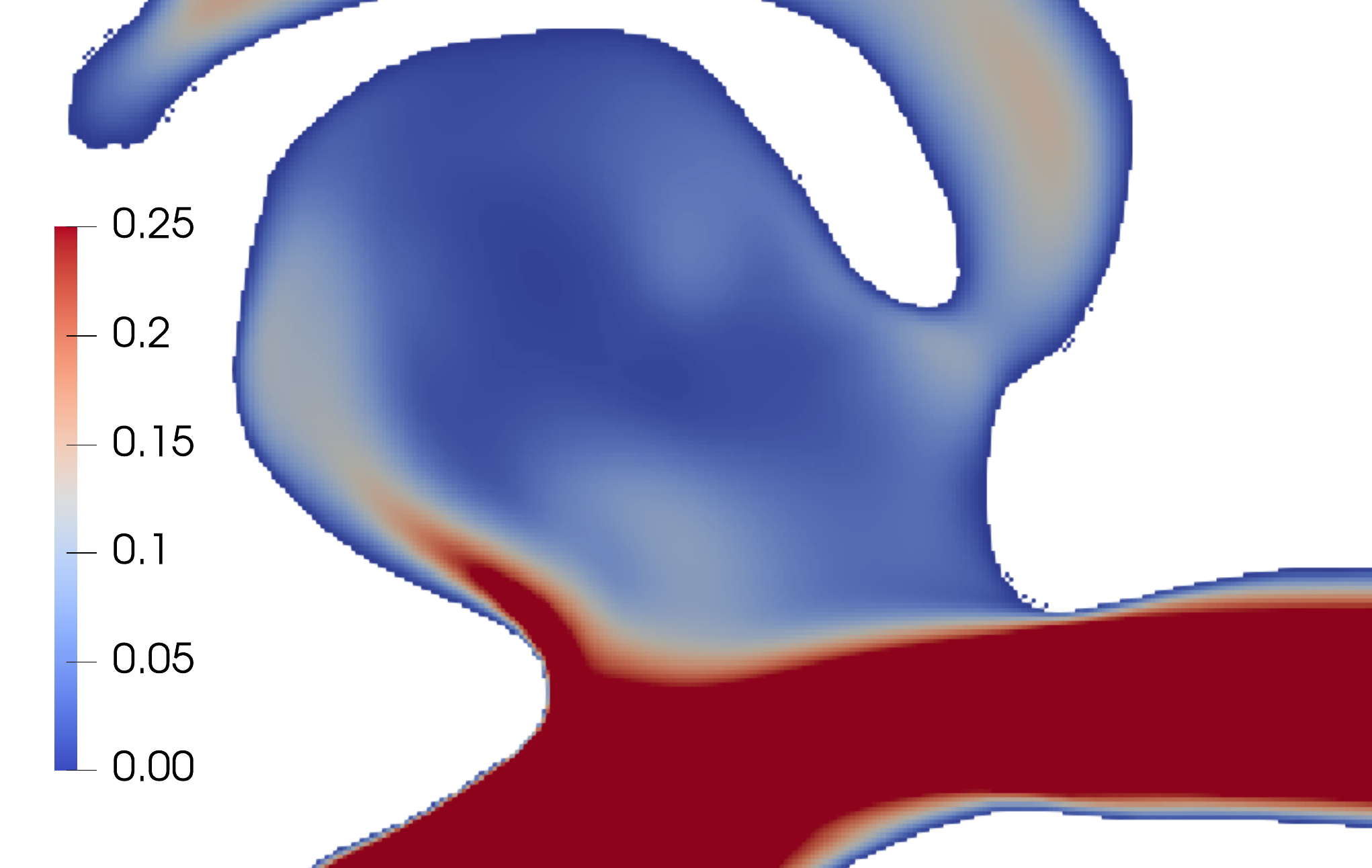}}%
    \end{minipage}
    \caption{Velocity magnitude $\|\vec{u}\|$ [m/s] during systole ($t \approx 1.23$\;s) in the aneurysm with fully resolved coil of 20\% packing density \emph{(left)} and with averaged porosity field \emph{(right)}.}
    \label{fig:vanse20}
\end{figure}
\begin{figure}[p]
    \begin{minipage}{\textwidth}
        \leftfigure{\includegraphics[width=.47\textwidth, viewport=370 350 1350 970, clip=true]{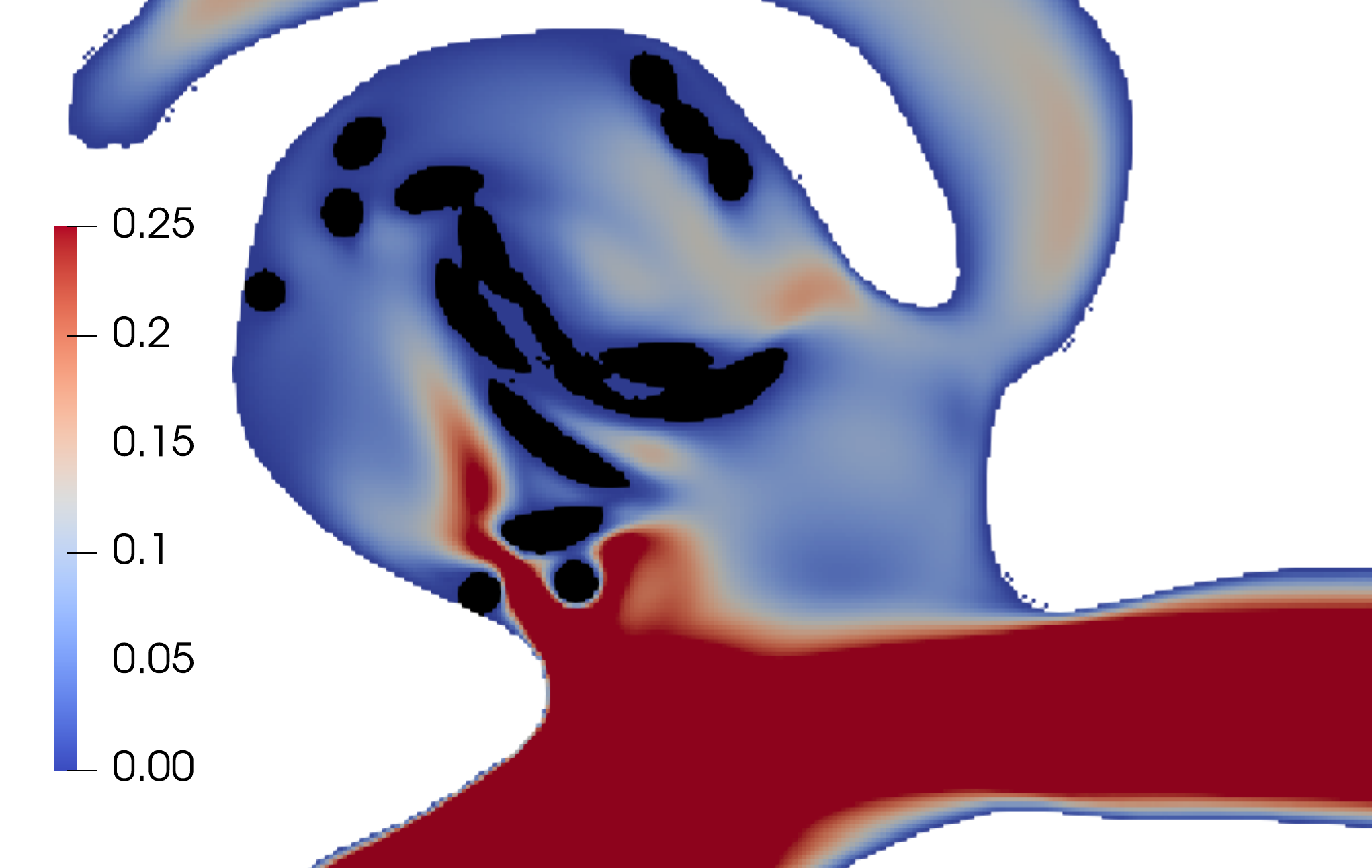}}%
        \hfill%
        \rightfigure{\includegraphics[width=.47\textwidth, viewport=370 350 1350 970, clip=true]{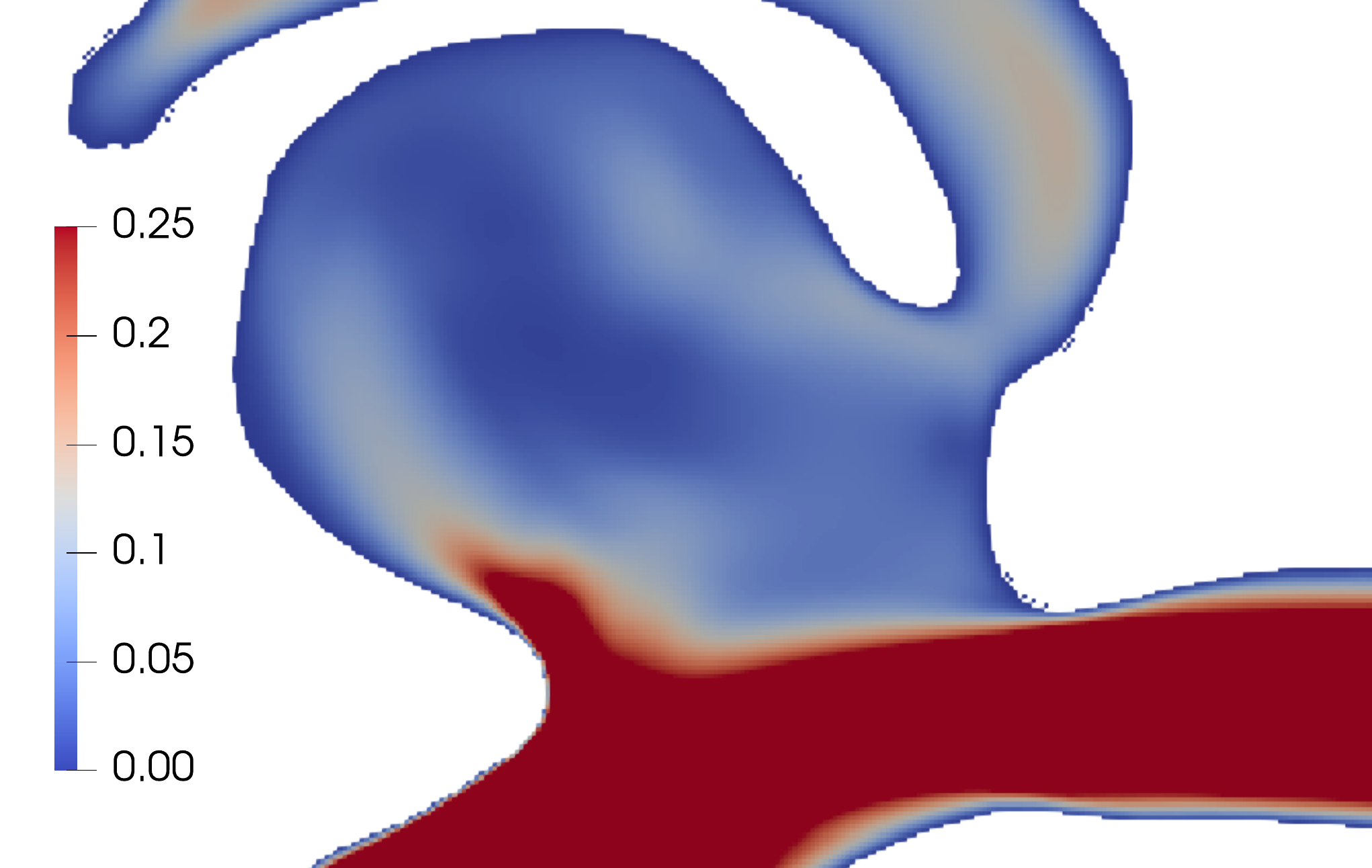}}%
    \end{minipage}
    \caption{Velocity magnitude $\|\vec{u}\|$ [m/s] during systole ($t \approx 1.23$\;s) in the aneurysm with fully resolved coil of 25\% packing density \emph{(left)} and with averaged porosity field \emph{(right)}.}
    \label{fig:vanse25}
\end{figure}

In a first step, we compare the fully resolved flow simulations in the aneurysm without and with treatment by coiling, see \cref{fig:none,fig:vanse15,fig:vanse20,fig:vanse25} (left subfigures).
In all fully resolved simulations, the majority of the fluid flows along the main (bottom) vessel.
The remaining flow enters the aneurysm without coil in an unimpeded manner, whereas the coil significantly reduces the flow in the aneurysm, while still allowing for flow through the secondary (upper) vessel, thus indicating a significant improvement by the treatment.

For the VANSE simulations depicted in the right subfigures, the averaged porosity field leads to good agreement with the fully resolved simulations.
In particular, only local fine-scale effects are lost due to the averaging, but the overall flow behavior is the same.
As expected, the flow reduction is more pronounced with increasing packing density.
As packing density varies in these examples from 15\% to 25\%, this indicates the validity of the model over the whole range encountered in practice.

\runinhead{Wall Shear Stress Reduction.}

Besides the flow reduction, a main quantity of interest associated with the formation and rupture of aneurysms is the wall shear stress $\vec{\tau} = \tens{\sigma} \vec{n} - (\vec{n} \cdot \tens{\sigma} \vec{n}) \vec{n}$, where $\vec{n}$ \update{is} the outer normal of the vessel wall.
This is computed based on the deviatoric stress $\tens{\sigma}$ in the flow domain by extrapolation to the wall mesh and subsequent projection onto the local mesh-tangential plane \cite{Krueger2017, Matyka2013}.

\begin{figure}[tb]
    \begin{minipage}{\textwidth}
        \leftfigure{\includegraphics[width=54mm, viewport=60 180 1100 930, clip=true]{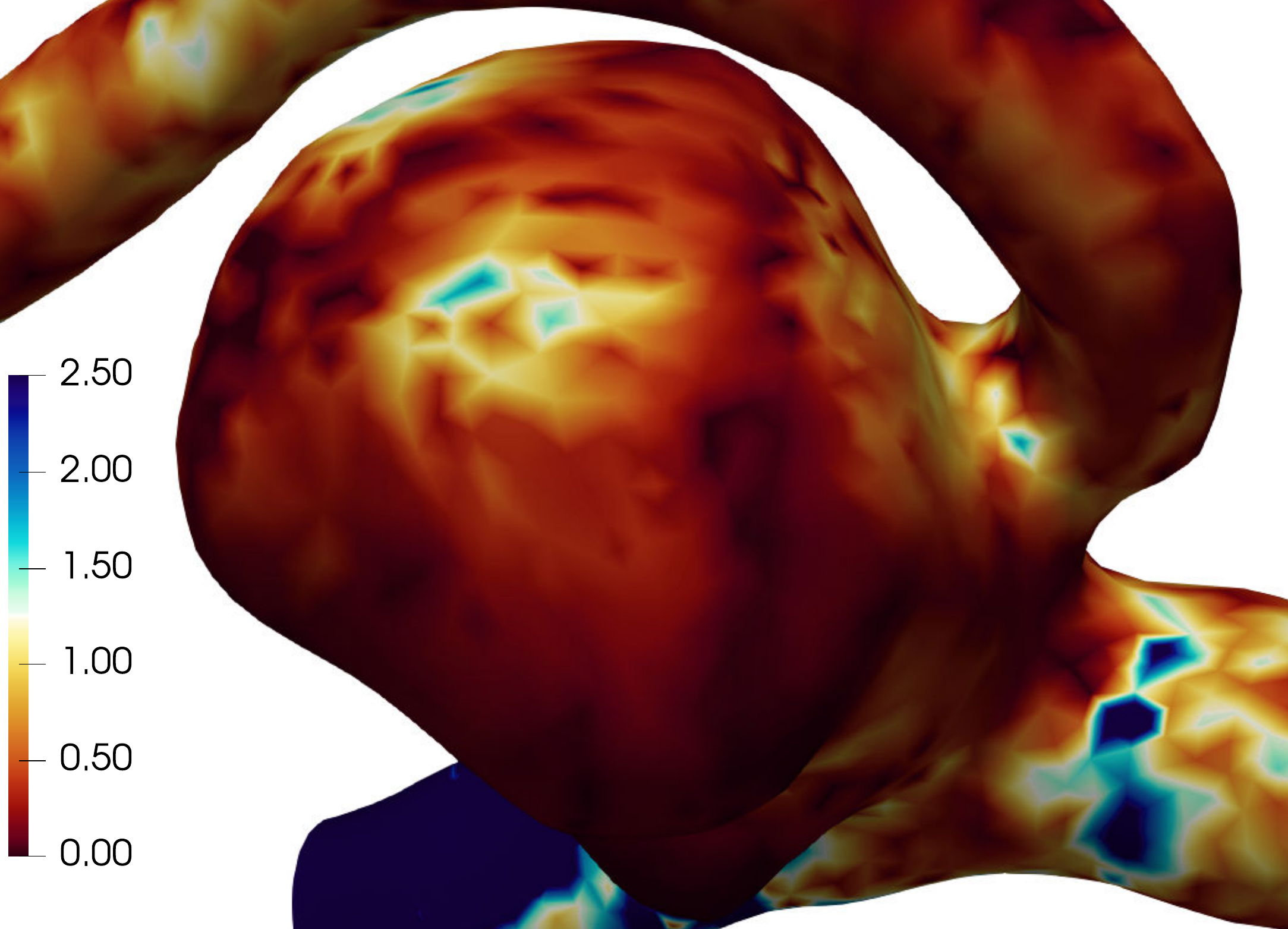}}%
        \hfill%
        \rightfigure{\includegraphics[width=54mm, viewport=60 180 1100 930, clip=true]{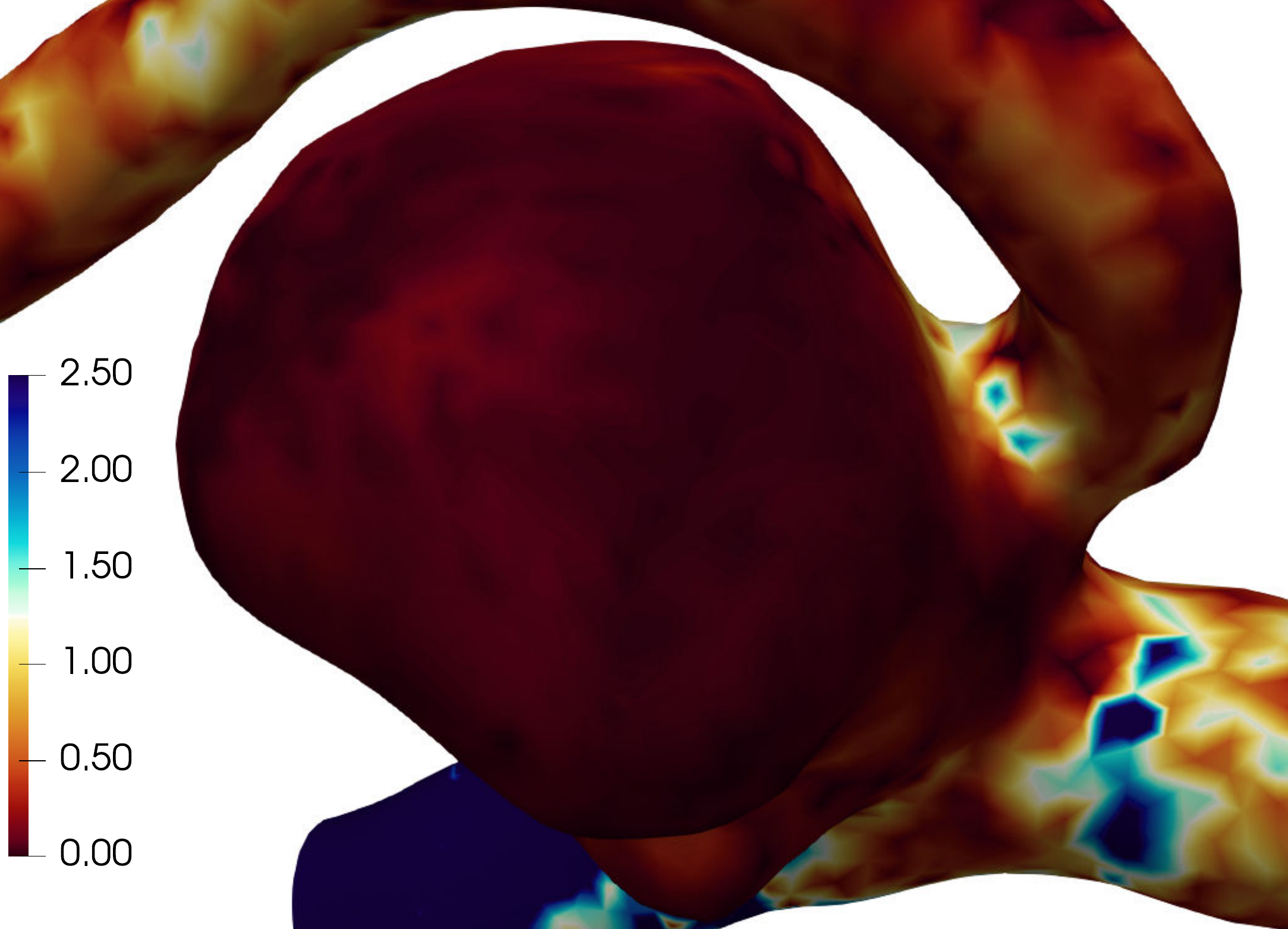}}%
    \end{minipage}
    \caption{Wall shear stress magnitude $\|\vec{\tau}\|$ [Pa] during systole ($t \approx 1.23$\;s) in the aneurysm without coil \emph{(left)} and with averaged coil of 25\% packing density \emph{(right)}.}
    \label{fig:3Da}
\end{figure}

The results during systole depicted in \cref{fig:3Da} show that the wall shear stress is drastically reduced by the treatment.
While the wall shear stress in the untreated aneurysm reaches a magnitude of about $2$\;Pa, all coiled aneurysms show a magnitude below $1.5$\;Pa and hence a lower risk of rupture.
Additionally, the wall shear stress magnitude differs only slightly for the three packing densities depending on the actual coil placement.

\begin{figure}[tb]
    \begin{minipage}{\textwidth}
        \leftfigure{\includegraphics[width=54mm]{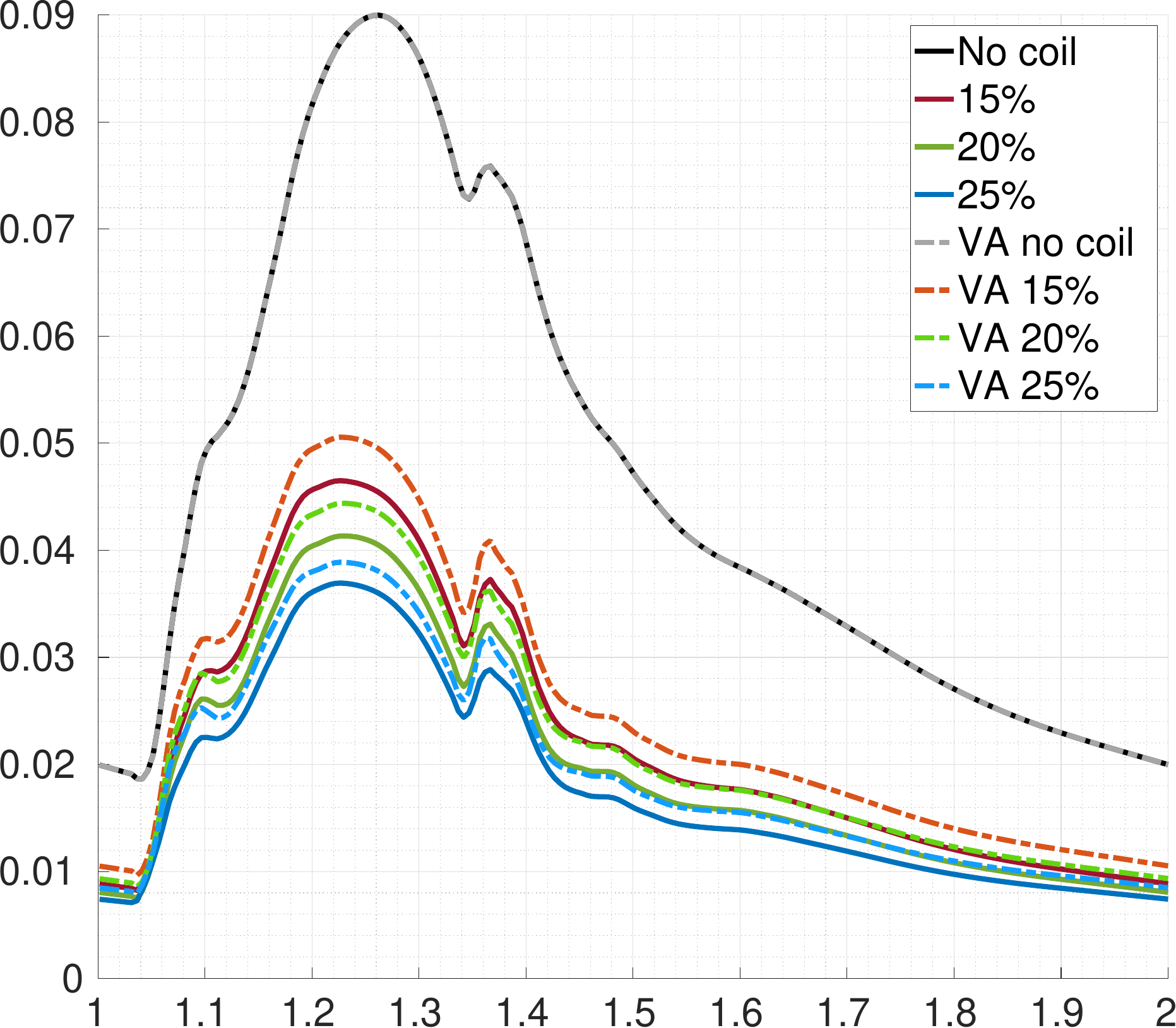}}%
        \hfill%
        \rightfigure{\includegraphics[width=54mm]{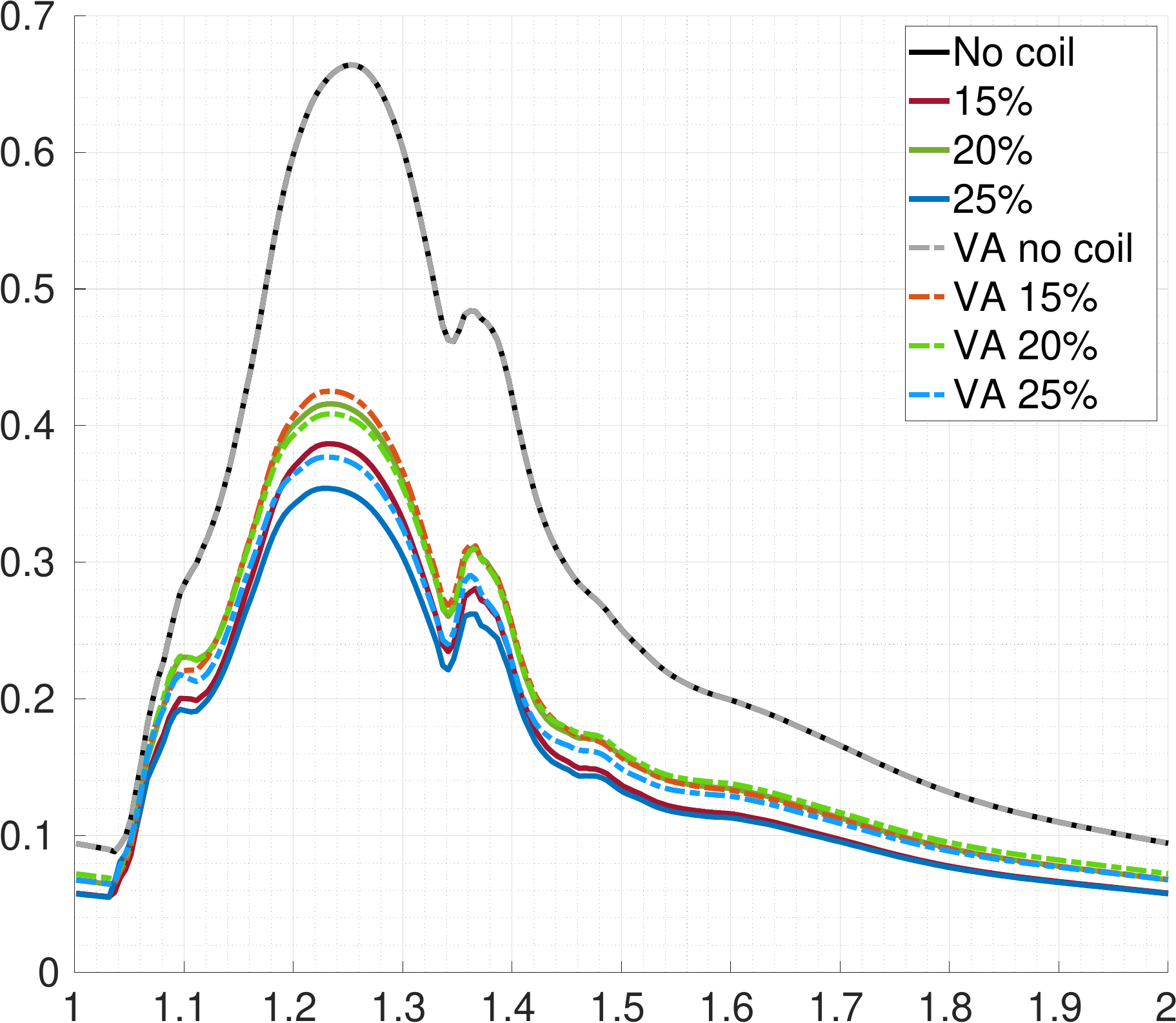}}%
    \end{minipage}
    \smallskip
    \caption{Average \update{in the aneurysm of the} velocity magnitude $\|\vec{u}\|$ [m/s] \emph{(left)} and of the wall shear stress magnitude $\|\vec{\tau}\|$ [Pa] \emph{(right)} plotted over the time $t$ [s] for the fully resolved simulations without and with coils of different packing density, and for the volume-averaged (VA) ones.}
    \label{fig:comparison}
\end{figure}

Furthermore, we compute the average \update{of the magnitudes of velocity $\|\vec{u}\|$ and wall shear stress $\|\vec{\tau}\|$ over the cut-off aneurysm region (colored light green in \cref{fig:none}b) and its intersection with the boundary of the original geometry, respectively}.
The time series over the second heart beat is depicted in \cref{fig:comparison} for the fully resolved and volume-averaged coils.
The flow magnitude is reduced by about 50\% and the wall shear stress magnitude by about 40\%.
Increasing the packing density of the coil reduces the overall flow in the aneurysm, but seems to \update{only slightly affect} the wall shear stress.
This might as well be due to the different coil distribution.
\medskip

In summary, this workflow allows for patient\update{-}specific assessment of the flow due to potential treatment.
In view of the large uncertainties \update{in} the exact coil placement as well as the required \update{fast} decision for or against treatment, accounting for an anisotropic permeability seems rather unnecessary in clinical practice, but might further increase the accuracy at the price of much slower simulations.
An analysis of this extension as well as a detailed parameter calibration and uncertainty quantification remain open question to be addressed in future research.

\begin{acknowledgement}
We thank Julian Schwarting and Jan Kirschke at the Neuro-Kopf-Zentrum of the Klinikum rechts der Isar (TU Munich) for providing the CT scan of the aneurysm.
\end{acknowledgement}

\ethics{Competing Interests}{
    The authors gratefully acknowledge the financial support provided by the German Science Foundation (DFG) under project number 465242983 within the priority programme ``SPP 2311: Robust coupling of continuum-biomechanical in silico models to establish active biological system models for later use in clinical applications -- Co-design of modeling, numerics and usability'' (WO 671/20-1) and (WO 671/11-1).
}

\bibliographystyle{spmpsci}
\bibliography{references}
\end{document}